\documentclass[11pt]{article}
\usepackage{amsmath,amsthm,amsfonts}
\newcommand{\Ga}{\Gamma}
\newcommand{\la}{\lambda}
\newcommand{\eps}{\varepsilon}
\newcommand{\C}{\mathbb C}
\newcommand{\Z}{\mathbb Z}
\newcommand{\N}{\mathbb N}
\newcommand{\Y}{\mathbb Y}

\newcommand{\Ind}{\operatorname{Ind}}
\newcommand{\Cent}{\operatorname{Cent}}
\newcommand{\Char}{\operatorname{Char}}
\newcommand{\mmod}{\operatorname{mod}}

\newcommand{\dfn}{\it}

\newtheorem{xprop}{Proposition}[section]
\newtheorem{xthm}[xprop]{Theorem}
\newtheorem{xlem}[xprop]{Lemma}
\newtheorem{xcor}[xprop]{Corollory}
\theoremstyle{remark}
\newtheorem{xrem}[xprop]{Remark}
\newtheorem{xdf}[xprop]{Definition}

\newtheorem{xexample}[xprop]{Exapmle}

\begin{document}
\begin{center}
{\Large \bf Traces on Infinite-Dimensional Brauer Algebras%
\footnote{Partially supported by the SRDF grant RUM1-2622-ST-04 and INTAS grant 03-51-5018.}.}\\[0.5in]
\vspace{2pt}%
{\large\bf A.~Vershik, P.~Nikitin}
\\[6pt]
St.Petersburg Department of Steklov Institute of Mathematics, \\
27 Fontanka, St.Petersburg 191023, Russia
\\[4pt]
E-mail: {vershik@pdmi.ras.ru, pnikitin0103@yahoo.co.uk}

\end{center}

\begin{abstract}
We describe the central measures for the random walk on graded graphs. Using this
description, we obtain the list of all finite traces on three infinite-dimensional
algebras: on the Brauer algebra, on the partition algebra, and on the walled Brauer
algebra. For the first two algebras, these lists coincide with the list of all finite
traces of the infinite symmetric group. For the walled Brauer algebra, the list of finite
traces coincide with the list of finite traces of the square of the latter group.
\end{abstract}

\section{Brauer algebras and the pascalized graphs}
Consider the diagonal action of complex orthogonal group $O_k(\C)$ acts on tensor power $V^{\otimes
n}$ of the space $V={\Bbb C}^k$:
$$
M\cdot (v_1\otimes\ldots\otimes v_n) = Mv_1\otimes\ldots\otimes
Mv_n, \quad M\in O_k(\C).
$$
R.~Brauer (see~\cite{Brauer, Weyl}) defined an family of finite dimensional algebra $Br_n(k)$ {\it
the Brauer algebra},depending on complex parameter $k$ and positive integer $n$. For integer $k
\geq k$ the algebra $Br_n(k)$ is isomorphic to centralizer of above described action of $O_k(\C)$.
For $k\in\ \C \mid \Z \cup\{n \dots\}$ and $n$ fixed these algebras are semisimple and pairwise
isomorphic \cite{Wenzl}. For now on we consider only these number of parameter $k$ and denote the
corresponding algebras by $Br_n$ omitting $k$ in the notation. For other $k \in \Z_- \cup \{1,\dots
n-1\}, k<n$ the algebra  $Br_n(k)$ is not semisimple.

We shall also study {\it the walled Brauer algebra} $Br_{n,m}(k)$, $n,m\in\Z_+$. The history of
this algebras is as follows. V.~Turaev (\cite{Turaev}), was first to define it by presentation, he
also pointed to the first author that its dimension is $(n+m)!$ and that this algebra resembles the
symmetric group. The walled Brauer algebra was defined independently in the work of K.~Koike
\cite{Koike}, and later it was studied in \cite{Benkart} as the centralizer of the diagonal action
of the group $GL_k(\C)$ on tensor space $V^{\otimes n}\bigotimes V^{*\otimes m}$. It was clear from
its diagrammatic definition that this algebra is subalgebra of the Brauer algebra. The walled
Brauer algebra is also semisimple and pairwise isomorphic for the generic $k$, $k\in\{x\in\C\mid
x\not\in\Z\} \cup \{x\in\Z \mid |x|\ge m+n\}$ (see \cite{Nikitin} for details). Here we again also
consider only these generic values and omit $k$ in the notation, $Br_{n,m}=Br_{n,m}(k)$.

P.~Martin introduced {\it the partition algebras} $Part_n(k)$, $n\in\Z_+$, $k\in\C$
(see~\cite{Martin}). Algebras $Part_{2n}(k)$, $Part_{2n+1}(k)$ for sufficiently large
$k\in\N$ are isomorphic to the centralizers of the diagonal action of the subgroups
$S_k\subset GL_k(\C)$, $S_{k-1}\subset GL_k(\C)$ on tensor space $V^{\otimes n}$. For the
generic $k\in\{x\in\C\mid x\not\in\Z\} \cup \{x\in\N\mid x\ge 2n-1\}$ and fixed $n$
partition algebras are semisimple and pairwise isomorphic, we will denote them by
$Part_n$.

Each finite-dimensional algebra concerned include the ideal $J$ (with an appropriate subscript)
spanned by all noninvertible standard generators of the corresponding algebra
(see~~\cite{Wenzl,Nikitin,Martin}) and the following holds:
\begin{equation}\label{eq:factors_finite_dim}
Br_n/J_n\cong \C[S_n],\quad Br_{n,m}/J_{n,m}\cong \C[S_n\times S_m],\quad
  Part_{2n}/J_{2n}\cong Part_{2n+1}/J_{2n+1} \cong \C[S_n].
\end{equation}

The algebras $Br_n$, $Br_{n,m}$, and $Part_n$ form the inductive families with natural embeddings.
This permit one to consider the inductive limits $Br_\infty=\varinjlim Br_n$,
$Br_{\infty,\infty}=\varinjlim Br_{n,m}$, and $Part_\infty=\varinjlim Part_n$, which are {\it the
locally semisimple (l.s.) algebras}. From the combinatorial point of view, each l.s. algebra (i.e.,
the inductive limit of the finite-dimensional semisimple algebras) $A=\varinjlim A_n$ is completely
determined by its {\it Bratteli diagram (branching graph)} $\Gamma(A)$. Recall that the {\it
Bratteli diagram} of the inductive family $\{A_n\}$ of finite-dimensional semisimple algebras is a
$\Z_+$-graded graph, constructed in the following way. Its vertices of the $n$-th level are
enumerated by the simple modules of the algebra $A_n$, and the edges between the $(n-1)$-th and
$n$-th levels correspond to the decomposition of the simple $A_n$-modules treated as
$A_{n-1}$-modules. (See surveys \cite{VershikKerov, Voiculescu} for the definitions concerning l.s.
algebras.) The problem of finding the Bratteli diagram of an l.s. algebra has much in common with
the problem of finding the spectrum of a commutative algebra.

Now we shall give the definition of the {\it paskalized graph}. It will be convenient for
us to use it for description of the Bratteli diagrams for the Brauer algebras and the
partition algebra. This definition includes the Jones basic construction, the main tool
in the finite-dimensional situation, see the papers of Jones and Wenzl \cite{Jones,
Wenzl, GdHJ}, the survey of Ram and Halverson \cite{HalversonRam} about the partition
algebras, and literature there in. Suppose $\Ga$ is a $\Z_+$-graded locally finite graph
with a single vertex on the $0$ level and without dangling vertices, and by $\Ga_k$
denote its set of vertices of the $k$-th level, $k\in \Z_+$. We also write $|\la|=i$ for
$\la\in\Ga_i$, and we write $\la\nearrow\nu$ ($\la\searrow\nu$) if the vertex $\nu$
follows the vertex $\la$ (precedes the vertex $\la$). Let us define the $\Z_+$-graded
graph $\Pi(\Ga)$. We set the $k$-th level $\Pi(\Ga)_k$ to be the union of the set $\Ga_k$
and of the sets $\Ga_i$ for all the previous levels of the same parity. We denote the
vertices in $\Pi(\Ga)_k$ by $(k,\la)$, where $\la\in \Ga_i$, $i\le k$, $k-i=0(\mmod 2)$.
We define the edges of the graph $\Pi(\Ga)$ as follows:
\begin{equation}\label{eq:edges}
(k,\la)\nearrow(k+1,\nu)\quad\Leftrightarrow\quad\la\nearrow\nu\text{ or }\la\searrow\nu.
\end{equation}
\begin{xdf}
The graph $\Pi(\Ga)$ is called the {\dfn pascalized graph} $\Ga$.
\end{xdf}
There is no difficulty to see that we can obtain the pascalized graph $\Pi(\Ga)$
iterating the following construction. Let us reflect the $(k-2)$-nd level of the
pascalized graph with respect to $(k-1)$-st level, add the set $\Ga_k$ to the result with
the edges that join the levels $\Ga_{k-1}$ and $\Ga_k$, and we get the $k$-th level of
the pascalized graph.

Obviously, a graph $\Ga$ is a subgraph of the graph $\Pi(\Ga)$, which constitutes a very
small part of the whole pascalized graph.

\begin{xexample}\label{examp:Pascal}
Consider the graph $\Ga^0$ with the set of vertices equal to $\Z_+$, and with edges
joining vertices $n$ and $(n+1)$, $n\in\Z_+$. Then $\Pi(\Ga^0)$ is the "half"\ of the
Pascal graph, $\Pi(\Ga^0)_{2k} = \{0,2,\dots,2k\}$, $\Pi(\Ga^0)_{2k+1} =
\{1,3,\dots,2k+1\}$, with the edges $(0,1)$, and $(i,i-1)$, $(i,i+1)$ for $i>0$; this
shows the origin of our definition. Note that the graph $\Pi(\Ga^0)$ correspond to the
Temperley-Lieb algebra (see~\cite{GdHJ}).
\end{xexample}

By $Y_\Ga$ we denote the set of all paths on a graph $\Ga$, in particular, we have
$Y_\Ga\subset Y_{\Pi(\Ga)}$. Space of all paths of a graph possesses the structure of a
totally disconnected compact set, the basis of open-closed sets consists of the cylinder
sets. We denote by $C_d=C^\Ga_d$ the cylinder that consists of all paths through the
vertex $d$.

Using \eqref{eq:edges}, we see that every path in the pascalized graph $\Pi(\Ga)$ is
uniquely determined by the sequence $\{x_n\}_{n=0}^\infty$ of the vertices of the graph
$\Ga$ such that each two neighbour vertices $x_n$, $x_{n+1}$ are the neighbours in the
graph $\Ga$, i.e., $x_n\nearrow x_{n+1}$ or $x_n\searrow x_{n+1}$. In other words, we can
say that a path in the pascalized graph $\Pi(\Ga)$ is a trajectory of the random walk on
the graph $\Ga$, therefore we can say that the graph $\Pi(\Ga)$ is the graph of the
random walk on the graph $\Ga$.

For example, the homogeneous tree of degree $2k$ with the marked vertice, which turnes
the tree into the $\Z_+$-graded graph, is the Caley graph of the free group with $k$
generators; and the corresponding pascalized graph is the graph of the random walk on the
free group.

We can consider the ideals $J$ for the infinite-dimensional algebras concerned in the
same manner as for the finite-dimensional case, and we get the similar picture:
\begin{equation*}
Br_\infty/J\cong \C[S_\infty],\quad Br_{\infty,\infty}/J\cong \C[S_\infty\times
S_\infty], \quad Part_\infty/J\cong \C[S_\infty].
\end{equation*}
We will see now that the Bratteli diagrams of the algebras $Br_\infty$,
$Br_{\infty,\infty}$, and $Part_\infty$ are the pascalized Bratteli diagrams for the
corresponding factoralgebras.

Recall, that the Bratteli diagram of the infinite symmetric group is the Young graph
$\Y$(see, for example, \cite{VershikKerov}). Now we will formulate in the convenient
terms known results about the algebras concerned. The Bratteli diagram for the family of
the Brauer algebras was obtained by H.~Wenzl (see~\cite{Wenzl}):
\begin{xthm}
The Bratteli diagram of the Brauer algebra $Br_\infty$ for the generic parameter is the
pascalized Young graph $\Y$, $\Ga(Br_\infty)=\Pi(\Y)$.
\end{xthm}

The Bratteli diagram for the walled Brauer algebra was obtained independently in
\cite{Leduc, Nikitin}. Consider the l.s. algebra $\C[S_\infty\times S_\infty]$ as the
inductive limit of the finite-dimensional algebras
$$
\C[S_{0}\times S_{0}]\subset\C[S_{1}\times S_{0}]\subset\C[S_{1}\times
S_{1}]\subset\C[S_{2}\times S_{1}]\subset\dots\subset\C[S_{[n+1/2]}\times
S_{[n/2]}]\subset\dots,
$$
and denote by $\bar{\Y}$ the corresponding Bratteli diagram, which can be simply
constructed using the graph $\Y$.
\begin{xthm}
The Bratteli diagram of the walled Brauer algebra $Br_{\infty,\infty}$ for the generic
parameter is the pascalized graph $\bar{\Y}$, $\Ga(Br_{\infty,\infty}) = \Pi(\bar{\Y})$.
\end{xthm}

P.~Martin described the Bratteli diagram for the family of the partition algebras
(see~\cite{Martin}). Consider the Young graph $\bar{\bar{\Y}}$ with each level repeated
twice, corresponding to the family of symmetric groups
$$
G_0\subset G_1\subset G_2\subset\dots,\qquad G_{2i}=G_{2i+1}=\C[S_i],\quad i\in\Z_+.
$$
\begin{xthm}
The Bratteli diagram of the partition algebra $Part_\infty$ for the generic parameter is
the pascalized graph $\bar{\bar{\Y}}$, $\Ga(Part_\infty)=\Pi(\bar{\bar{\Y}})$.
\end{xthm}

\section{Central measures on the pascalized graphs}

One of the main questions in the theory of l.s. algebras is to describe finite traces and
the $K$-functor of an algebra. A linear functional $f:A\to\C$ is called a {\it (finite)
trace} on a *-algebra $A$ if the following conditions hold
\begin{enumerate}
 \item $f(1)=1$;
 \item $f(ab)=f(ba)$,\quad $a,b\in A$;
 \item $f(a^*a)\ge0$,\qquad $a\in A$.
\end{enumerate}
By $\Char(A)=\Char(\Ga(A))$ denote the set of all finite traces on $A$. A trace is called
{\it indecomposable} if it is not equal to a nontrivial convex sum of traces.

A probability Borel measure on space $Y_\Ga$ of paths of a $\Z_+$-graded infinite graph
$\Ga$ is called a {\it central measure} if it induces the uniform conditional measures on
finite sets of paths that coincide from some place. Denote the set of all central
measures on space of paths of a graph $\Ga$ by $\Cent(\Ga)=\Cent(Y_\Ga)$. A central
measure is called {\it ergodic} if it is not equal to a nontrivial convex sum of central
measures. There is a well-known correspondence between the set of finite traces on an
l.s. algebra $A$ and the set of central measures on space of paths of the Bratteli
diagram $\Ga(A)$:
$$
\Char(A)\leftrightarrow \Cent(\Ga(A))
$$
(see, for example, \cite{VershikKerov, Voiculescu}). Indecomposable traces here
correspond to the ergodic central measures.

We will reduce the study of traces on a pascalized graph $\Pi(\Ga)$ (for some algebras)
to the study of traces on an initial graph $\Ga$, and then we will use this result to
find the list of traces on the Brauer algebras and on the partition algebra. {\it We
prove that, under certain conditions on the Bratteli diagram, every central measure on
the pascalized graph is nonzero only on the initial graph as a subgraph of the pascalized
one. In other words, every central measure on the pascalized graph coincide with a
central measure on the initial graph.} In particular, we prove that there is a one-to-one
correspondence between the traces on algebras $Br_\infty$ and $Part_\infty$ and the
traces on the infinite symmetric group $S_\infty$, and that the traces on
$Br_{\infty,\infty}$ are in one-to-one correspondence with the traces on the group
$S_\infty\times S_\infty$.

The list of central measures for the infinite-dimensional Brauer algebra was given
without proof in the paper \cite{Kerov2} by S.~Kerov; traces on the infinite-dimensional
walled Brauer algebra and on the partition algebra, probably, were considered nowhere.
Also there was no study of the set of the full measure for the central measures
concerned.

Let $\Ga$ be a $\Z_+$-graded graph. By $\dim(d;d')$ denote the number of paths from the
vertex $d$ to the vertex $d'$ in $\Ga$, by $\dim(d)$ denote the number of paths from the
initial vertex of the graph to the vertex $d$. Recall the {\it ergodic method} for
finding the central measures:
\begin{xthm}[\cite{VershikKerov}]\label{thm:center_measure}
Consider a $\Z_+$-graded graph, an ergodic central measure $\mu$ on it, and the set $S$
of paths of the form $s = (s_0\nearrow s_1\nearrow\dots\nearrow s_n\nearrow\dots)$ such
that for every vertex $d$ holds
\begin{equation}\label{eq:center_measure}
\mu(C_d) = \lim_{n\to\infty}\frac{\dim(d)\cdot\dim(d;s_n)}{\dim s_n};
\end{equation}
then $\mu(S) = 1$.
\end{xthm}

Hence, for finding the central measures it is sufficient to describe all the limits
\eqref{eq:center_measure}. Now we use this method for pascalized graphs.

\begin{xlem}\label{lem:mu(2,emptyset)}
If the following condition hold for the vertices of a graph $\Pi(\Ga)$:
\begin{equation}\label{eq:lim_frac}
  \lim_{n\to\infty}\max_{\la,\ |\la|<n}\frac{\dim(n-2,\la)}{\dim(n,\la)} = 0,
\end{equation}
then for every ergodic central measure $\mu$ on $\Pi(\Ga)$ and for every vertex
$(n_0,\la)\in\Pi(\Ga)_{n_0}$, $|\la|<n_0$, holds $\mu(C_{(n_0,\la)}) = 0$.
\end{xlem}
\begin{proof}
If we take an ergodic central measure $\mu$ and a path $s$ from
Theorem~\ref{thm:center_measure}, then the following holds
$$
\mu(C_{(n_0,\la)}) = \lim_{n\to\infty}\frac{\dim(n_0,\la)\cdot\dim((n_0,\la);(n,s_n))}
{\dim(n,s_n)}.
$$
Consider a vertex $(n,s_n)$ on the path from the vertex $(n_0,\la)$, $|\la|<n_0$. It is
obviously true, that $|s_n|<n$. Using the parity, we get $|\la|\le n_0-2$, $|s_n|\le
n-2$, therefore the vertices $(n_0-2,\la)$ and $(n-2, s_n)$ are well defined. There is a
bijection between the paths from the vertex $(n_0,\la)$ to $(n,s_n)$ and the paths from
$(n_0-2,\la)$ to $(n-2,s_n)$, we simply increase by $2$ the number of the level for every
vertex in the path. Thus
\begin{equation*}
 \mu(C_{(n_0,\la)}) =
    \dim(n_0,\la)\cdot\lim_{n\to\infty}\frac{\dim((n_0-2,\la);(n-2,s_n))}{\dim(n,s_n)}\le
    \dim(n_0,\la)\cdot\lim_{n\to\infty}\frac{\dim(n-2,s_n)}{\dim(n,s_n)},
\end{equation*}
and the lemma is proved.
\end{proof}

Consider an inductive family of finite-dimensional semisimple algebras
$$
A_0\cong\C\subset A_1\subset\dots\subset A_l\subset\dots,
$$
and set $a_{l+1}=[\dim A_{l+1}/\dim A_l]$, $l\in\Z_+$, where $\dim A_l$ is the dimension
of the algebra $A_l$.
\begin{xlem}\label{lem:multiplicity}
Let $\Ga$ be the branching graph of the family
$$
A_0\cong\C\subset A_1\subset\dots\subset A_l\subset\dots.
$$
If for every $l\in\Z_+$, and every simple module $\la\in \Ga_l$ holds
$$
\dim(\Ind_{A_l}^{A_{l+1}}\la)=a_{l+1}\dim\la,
$$
then the ratio of the dimensions $\dim(n,\la)/\dim\la$ depend only on the number $n$ of
the level of the vertex in the pascalized graph and on the number $|\la|$ of the level of
the vertex in the initial graph. In other words, there exist numbers $M(n,l)$ for $l\le
n$, $n-l=0(\mmod 2)$ such that
$$
\dim(n,\la)=M(n,|\la|)\dim\la
$$
for every simple module $(n, \la)\in \Pi(\Ga)_n$. The numbers $M(n,l)$ can be defined
inductively as follows:
\begin{align}\label{eq:def_M}
M(n,n) &= 1;\quad M(2n+2,0)=M(2n+1,1),  \quad n=0,1,2,\dots;\\
M(n,l) &= M(n-1,l-1) + a_{l+1} M(n-1,l+1),  \quad 0<l<n,\quad n=0,1,2,\dots\ .\notag
\end{align}
\end{xlem}
\begin{proof}
Using the branching rule in the graph $\Pi(\Ga)$ we get
$$
\dim(n,\la)=\sum_{\eta\nearrow\la}\dim(n-1,\eta)+\sum_{\nu\searrow\la}\dim(n-1,\nu).
$$
By the inductive assumption,
$$
\dim(n,\la)=M(n-1,|\la|-1)\sum_{\eta\nearrow\la}\dim\eta+M(n-1,|\la|+1)\sum_{\nu\searrow\la}\dim\nu.
$$
Equality $\sum_{\eta\nearrow\la}\dim\eta = \dim\la$ holds for every branching graph. In
our case $\sum_{\nu\searrow\la}\dim\nu
=\dim(\Ind_{A_{|\la|}}^{A_{|\la|+1}}\la)=a_{|\la|+1}\dim\la$, hence
$$
\dim(n,\la)=\biggl(M(n-1,|\la|-1)+ a_{|\la|+1}M(n-1,|\la|+1)\biggr)\dim\la,
$$
and the formula is proved. The case $|\la|=0$ can be considered in the same manner, the
case
$|\la|=n$ is obvious.%
\end{proof}
\begin{xrem}
The converse to Lemma~\ref{lem:multiplicity} can be proved in the same way.
\end{xrem}
\begin{xrem}
Let $\{G_n\}$ be the inductive family of finite groups, then the conditions of
Lemma~\ref{lem:multiplicity} hold for the group algebra $\C[G]=\cup\C[G_n]$, in this case
$a_{l+1}=[G_{l+1}:G_{l}]$. Example~\ref{examp:Pascal} above shows, that in general
situation the list of traces on a graph $\Ga$ and the list of traces on a pascalized
graph $\Pi(\Ga)$ can be quite different: there is only one trace on $\Ga^0$, and the
traces on the pascalized graph $\Pi(\Ga)$ are parametrized by the segment $[0,1/2]$.
\end{xrem}
\begin{xlem}
If $n\ge2$, then
\begin{multline*}
 \frac{M(2n-2,0)}{M(2n,0)}=\frac{M(2n-3,1)}{M(2n-1,1)}>\frac{M(2n-2,2)}{M(2n,2)}>\frac{M(2n-3,3)}{M(2n-1,3)}>\frac{M(2n-2,4)}{M(2n,4)}>\dots, \\
 \frac{M(2n-2,0)}{M(2n,0)}>\frac{M(2n-1,1)}{M(2n+1,1)}>\frac{M(2n-2,2)}{M(2n,2)}>\frac{M(2n-1,3)}{M(2n+1,3)}>\frac{M(2n-2,4)}{M(2n,4)}>\dots\ .
\end{multline*}
\end{xlem}
\begin{proof}
The equality $\frac{M(2n-2,0)}{M(2n,0)}=\frac{M(2n-3,1)}{M(2n-1,1)}$ follows from
\eqref{eq:def_M}.

Recall the following simple fact:
$$
\frac{A}{B}>\frac{C}{D};\
A,B,C,D,x>0\Rightarrow\frac{A}{B}>\frac{A+xC}{B+xD}>\frac{C}{D}.
$$
Using \eqref{eq:def_M}, we get $M(2n-2,2) = M(2n-3,1)+a_3M(2n-3,3)$, $M(2n,2) =
M(2n-1,1)+a_3M(2n-1,3)$. Combining this two facts, we obtain inequalities
$\frac{M(2n-3,1)}{M(2n-1,1)}>\frac{M(2n-2,2)}{M(2n,2)}>\frac{M(2n-3,3)}{M(2n-1,3)}$. The
rest of inequalities can be obtained in the similar manner.
\end{proof}
\begin{xcor}
Set $\eps(n) = 0$ for $n$ even and $\eps(n) = 1$ for $n$ odd. Then under the conditions
of Lemma~\ref{lem:multiplicity}, we have
$$
\max_{\la,\ |\la|<n}\frac{\dim(n-2,\la)}{\dim(n,\la)} =
\frac{M(n-2+\eps(n),0)}{M(n+\eps(n),0)}.
$$
\end{xcor}
\begin{xcor}
The sequence $\{\frac{M(2n,0)}{M(2n+2,0)}\}_{n=0}^\infty$ is decreasing.
\end{xcor}
So, under the conditions of Lemma~\ref{lem:multiplicity}, the limit \eqref{eq:lim_frac}
tends to zero iff the sequence $\{\frac{M(2n,0)}{M(2n+2,0)}\}_{n=0}^\infty$ tends to
zero.
\begin{xlem}
\begin{enumerate}
  \item $\lim\limits_{n\to\infty}\frac{M(2n,0)}{M(2n+2,0)}>0\quad\Leftrightarrow\quad\exists M, C
  \ \forall n\ M(2n,0)<C\cdot M^n$;
  \item $\lim\limits_{n\to\infty}\frac{M(2n,0)}{M(2n+2,0)} = 0\quad\Leftrightarrow\quad\forall M
  \ \exists C\ \forall n\ M(2n,0)>C\cdot M^n$.
\end{enumerate}
\end{xlem}
\begin{proof}
The decreasing sequence $\{m_n=\frac{M(2n,0)}{M(2n+2,0)}>0\}_{n=0}^\infty$ always tends
to some limit $m$, $\lim_{n\to\infty}m_n=m\ge 0$. If $m>0$, then
$M(2n,0)=\prod_{k=0}^{n-1}\frac{M(2k+2,0)}{M(2k,0)}<C(\frac 1 m)^n$. If $m=0$, then for
every $M>0$ there exists $N$ such that $\frac{M(2k,0)}{M(2k+2,0)}<\frac 1 M$ for $k>N$.
In the latter case $M(2n,0)=\prod_{k=0}^{n-1}\frac{M(2k+2,0)}{M(2k,0)}>C_1\cdot
M^{n-N}=C_2\cdot M^n$.
\end{proof}
\begin{xthm}\label{thm:main}
Under the conditions of Lemma~\ref{lem:multiplicity}
\begin{equation*}
  \lim_{n\to\infty}\max_{\la,\ |\la|<n}\frac{\dim(n-2,\la)}{\dim(n,\la)} = 0
  \quad\Leftrightarrow\quad   \sup_l\{a_l=[\dim A_{l}/\dim A_{l-1}]\}=\infty.
\end{equation*}
\end{xthm}
\begin{proof}
Suppose $\sup_l\{a_l\} = a$; then using \eqref{eq:def_M} one can simply check by
induction that $M(n,l)< (a+1)^n$ for every $l$.

If $\sup_l\{a_l\} = \infty$, then for every $M$ there exists an index $L$ such that
$a_L>M$. Suppose $2n>2L$ and consider the path
\begin{multline*}
(0,0)\nearrow(1,1)\nearrow\dots\nearrow(L,L)\nearrow(L+1,L-1)\nearrow(L+2,L)\nearrow(L+3,L-1)\nearrow\dots \\
\dots\nearrow(2n-L,L)\nearrow(2n-L+1,L-1)\nearrow(2n-L+2,L-2)\nearrow\dots\nearrow(2n,0).
\end{multline*}
Using induction, we obtain $M(L+2i,L)>M^i$, hence $M(2n,0)>M^{n-2L}=C M^n$.
\end{proof}

Thus we have described the conditions when the central measures does not change if we
pascalize the graph. Under the conditions of Theorem~\ref{thm:main}, we have
\begin{enumerate}
  \item every central measure on the pascalized graph $\Pi(\Ga)$ is nonzero only on the
  subgraph $\Ga\subset\Pi(\Ga)$, hence every such central measure coincides with a central measure on $\Ga$;
  \item consequently, there exists a bijection between the central measures $\Cent(\Ga)\leftrightarrow\Cent(\Pi(\Ga))$
  and, therefore, a bijection between the traces on the corresponding algebras.
\end{enumerate}

For the family of the Brauer algebras we have $a_l=\dim \C[S_l]/\dim \C[S_{l-1}] = l$,
for the walled Brauer algebras $a_l = [(l+1)/2]$, and for the partition algebras
$a_{2l}=l$, $a_{2l+1}=1$, this leads us to the following theorem.

\begin{xthm}\label{thm:main,measures}
Bratteli diagram for each of the algebras $Br_\infty$, $Br_{\infty,\infty}$,
$Part_\infty$ is the pascalized graph $\Pi(\Ga)$ for some graph $\Ga$. Central measures
on the graph $\Pi(\Ga)$ are nonzero only on the subgraph $\Ga\subset\Pi(\Ga)$, and there
is a bijection between these central measures and the central measures on $\Ga$.
\end{xthm}
\begin{xcor}
Every trace on each of the algebras $Br_\infty$, $Br_{\infty,\infty}$ and $Part_\infty$
is the lift of some trace on the corresponding factoralgebra $\C[S_\infty]$,
$\C[S_\infty\times S_\infty]$ и $\C[S_\infty]$.
\end{xcor}
Our description of the traces remain true for the deformation of the Brauer algebra ({\it
Birman-Wenzl algebra} and for the deformation of the walled Brauer algebra for the
generic parameters (about the deformations see \cite{BirmanWenzl, Murakami, Kerov2,
Turaev, KosudaMurakami, Halverson, Leduc}).

\section{About the $K_0$-functor}
Recall that the {\it $K_0(A)$-functor (Grothendick group)} of a locally semisimple
algebra $A$ is the inductive limit of the finitely generated infinite cyclic groups, the
generators of these groups are the vertices of the Bratteli diagram of the algebra $A$,
and the embeddings are also defined by this diagram. To see this description, we are to
identify the vertices of the branching graph and the corresponding simple
finite-dimensional modules of the finite-dimensional subalgebras. Then we consider the
induced modules, and thus we get the modules over the whole algebra $A$
(see~\cite{VershikKerov}). The vertex is called {\dfn infinitesimal}, if the cylinder of
the paths through this vertex is of measure zero for every central measure on the graph.
A projective module (and, in particular, a simple module) is called {\dfn infinitesimal},
if the corresponding vertex is infinitesimal, in other words, if all the finite traces
are zero on it. Hence we define the subgroup $I(A)=I$ of the infinitesimal modules in the
Grothendick group. One can easily see that simple infinitesimal modules generate $I$. For
the group algebra of the infinite symmetric group this subgroup is trivial. For the
algebras $Br_\infty$, $Br_{\infty,\infty}$ and $Part_\infty$
Theorem~\ref{thm:main,measures} says that the subgroup $I$ is not trivial: it is
generated by all the vertices in the difference $\Pi(\Ga)\backslash\Ga$ for the
corresponding graph $\Ga$. This leads to the following theorem:
\begin{xthm}
$$
K_0(Br_\infty)/I\cong K_0(\C[S_\infty]),\quad K_0(Br_{\infty,\infty})/I\cong
K_0(\C[S_\infty\times S_\infty]),\quad K_0(Part_\infty)/I\cong K_0(\C[S_\infty]).
$$
\end{xthm}
It seems that there are no works concerning the structure of the infinitesimal modules of
an l.s. algebra. These modules were studied in \cite{VershikKokhas} as the modules over
the complex group algebra of the group $SL(2,F)$ over the countable field $F$. Such
module can't be a submodule of the factor-representation with a trace (in particular, of
the regular representation), because a finite trace correspond to the finite-dimensional
representation or to the $II_1$-representation. This resembles the fact, that a
non-amenable group can possess an irreducible representation, which is not (weakly)
included in the regular representation.

\end{document}